\documentclass{article}[12pt] 

\usepackage{amssymb}
\usepackage[english]{babel}
\usepackage{amsfonts,epsfig}
\usepackage{epstopdf}
\addto\captionsenglish{}

\renewcommand{\i}{i}
\newcommand{\ptl}{\partial}

\newtheorem{problem}{Problem}

\DeclareMathAlphabet\mathbit
    \encodingdefault\rmdefault\bfdefault\itdefault
\DeclareOldFontCommand{\bi}{\normalfont\bfseries\itshape}{\mathbit}

\newcommand{\be}{\begin{equation}}
\newcommand{\ee}{\end{equation}}

\begin{document}

\title{Diffraction by an impedance strip I. Reducing diffraction problem to Riemann--Hilbert problems}

\author{A. V. Shanin, A. I. Korolkov}

\maketitle

\begin{abstract}
A 2D problem of acoustic wave scattering by a segment bearing impedance boundary conditions is considered.
In the current paper (the first part of a series of two) some preliminary steps are made, namely, the diffraction problem
is reduced to two matrix Riemann--Hilbert problems with exponential growth of unknown functions (for
the symmetrical part and for the antisymmetrical part). For this, the Wiener--Hopf problems are formulated, they
are reduced to auxiliary functional problems by applying the embedding formula, and finally the Riemann--Hilbert problems
are formulated by applying the Hurd's method.

In the second part the Riemann--Hilbert problems will be solved by a novel method of OE--equation.
\end{abstract}

\section{Introduction}

We study a 2D problem of diffraction by a segment bearing im\-pe\-dan\-ce boundary conditions
on both sides. This problem can be considered as a cross-section of a 3D problem of diffraction by an infinitely long strip with finite width and zero thickness.
The governing equation is the Helmhotz one, so the stationary problem is studied. No restriction is
imposed on the relation between the wavelength and the width of the strip (length of the segment).
The impedances of the sides are assumed to be equal.

The problem of diffraction by a segment has been studied extensively, but the vast majority of papers
is related to the case of ideal (Dirichlet or Neumann) boundary conditions. A problem with ideal boundary conditions (ideal segment) admits an application of separation of variables method in the elliptical coordinates.
As the result, the solution becomes expressed in terms of Matheu functions \cite{Sieger}.
However this solution seems not attractive for applications and for analytical studies.
Numerous attempts have been made to obtain a solution analogous to the Sommerfeld's formula for the
half-plane \cite{Somm}. A review of these attempts can be found in  \cite{Luneburg}.
Unfortunately, it has been found that the elegant approach of Riemann surface and Sommerfeld integral  cannot be successfully used for the segment problem.

A good practical way to treat the segment problem at least in the short-wave approximation is
the diffraction series approach. For an ideal segment this approach has been developed in
\cite{Sch,HasWein} and in many other papers.

Some mathematically important results for the ideal strip problem have been obtained in \cite{Williams,Latta,Gorenflo}.
The problem of diffraction by an ideal strip was reduced there to the inverse monodromy problem for a confluent Heun's equation. Thus the problem of diffraction by an ideal strip has been solved at least in the mathematical sense.
One of the authors contributed to this branch \cite{Shanin2001,Shanin2003a,Shanin2003b}.

The problem of diffraction by an impedance segment seems much more complicated. In the case of high frequencies the method of diffraction series can be applied to this problem \cite{Herman}. Otherwise one needs to solve an appropriate integral equation \cite{Senior} numerically. Also there exist some hybrid techniques, which combine both analytical and numerical approach. By using such techniques computational time may be significantly reduced  \cite{Burnside,Sahalos,Ikiz}. Besides, some approximate analytical methods, e.~g.\ an approximate Wiener--Hopf technique \cite{Serbest}
can be applied to this problem. Still the analytical theory of scattering by an impedance segment is far from being completed. Here we present some results that seem important and enable one to perform efficient calculations.

The first part of the paper describes the preliminary steps. Namely, the problem is formulated and symmetrized. After sym\-met\-ri\-za\-ti\-on, the sym\-met\-ri\-cal and the antisymmetrical problem are studied in parallel (they are slightly different). Following \cite{Nobl}, for each of these two diffraction problems a functional problem is formulated. Then, auxiliary functional problems are formulated. The {\em embedding formula\/} expressing the directivity in terms of the auxiliary solutions is derived. This embedding formula is useful since it represents the directivity (which is a function of the angle of incidence and the angle of scattering) as
a combination of functions depending on a single variable.

Method of embedding formula have been applied to many diffraction problems with different sets of auxiliary problems. In \cite{Williams} embedding formula was derived for diffraction by an ideal strip. Problems with grazing incidence were taken to generate auxiliary solutions. In \cite{Biggs1, Biggs2, Biggs3} embedding formula was obtained for diffraction by thin breakwaters using tricky manipulation with integral equations. Also embedding formula was derived  for planar cracks in \cite{Shanin2001b}. Edge Green's functions were used to generate auxiliary problems. In the current research we do not use this approach and just introduce auxiliary functional problems with a proper behaviour at infinity.

Then, following the procedure developed in \cite{Hurd} matrix Riemann--Hilbert problems are formulated for the auxiliary functional problems.

The second part of the paper will be dedicated to solving the matrix Riemann--Hilbert problems using a novel technique of the OE--equation.


\section{Formulation of diffraction problem}

Consider a 2D plane $(x,y)$. The scatterer is the segment $y = 0$, $-a< x < a$.
Everywhere outside this segment the Helmholtz equation is valid:
\begin{equation}
\Delta u + k_0^2 u = 0
\label{eq0101}
\end{equation}
where $u(x,y)$ is a field variable, and $k_0$ is a parameter.
We assume that $k_0$ has a vanishing positive imaginary part in order to use the limiting absorption principle.
The choice of time dependence is such that the wave traveling in the positive $x$-direction has the form
$e^{i k_0 x}$.

The total field is a sum of
the incident wave $u^{\rm in}$ and the scattered wave $u^{\rm sc}$:
\[
u = u^{\rm in} + u^{\rm sc},
\]
where
\begin{equation}
u^{\rm in} = \exp \{ - i k_0 (x \cos \theta^{\rm in} + y \sin \theta^{\rm in} ) \}
\label{eq0102}
\end{equation}
is a plane wave.
Here $\theta^{\rm in}$ is the angle of incidence; $0 \le \theta^{\rm in} \le \pi/2$.

The total field should be one-side continuous on the scatterer
and obey impedance boundary conditions on the faces of the scatterer:
\begin{equation}
\pm \frac{\ptl u}{\ptl y} (x , \pm 0) = \eta \, u (x , \pm 0) , \qquad -a < x < a.
\label{eq0103}
\end{equation}
Here $\eta$ is the impedance parameter. Energy conservation or dissipation condition requires
\begin{equation}
{\rm Im}[\eta] \le 0 .
\label{eq0104}
\end{equation}

The total field should obey Meixner's conditions near the vertices $(\pm a,0)$.
Namely, the integral of the ``energy''
combination $|\nabla u|^2 + |u|^2$ over any finite proximity of a vertex should be finite.
Later on, the Meixner's condition will be reformulated as a restriction imposed on the growth of the
field near the vertices.

The scattered field $u^{\rm sc}$ should also obey the Sommerfeld's radiation condition in the standard form:
\begin{equation}
\left(
\frac{\ptl u^{\rm sc}}{\ptl r} - i k_0 u^{\rm sc}
\right)
=
o(e^{i k_0 r} (k_0 r)^{-1/2}),
\label{eq0105}
\end{equation}
where
$r = \sqrt{x^2 + y^2}$.
Thus,
the scattered field for large $r$ can be written as follows:
\begin{equation}
u^{\rm sc} (r , \theta) = \frac{\exp\{i k_0 r\}}{\sqrt{2 \pi k_0 r}} S (\theta, \theta^{\rm in})
+ o(e^{i k_0 r} (k_0 r)^{-1/2}).
\label{eq0106}
\end{equation}
Here $\theta = \arctan (y/x)$, and $S(\theta, \theta^{\rm in})$ is the {\em directivity} of the
scattered field. This directivity should be found as the result of this research.


\section{Symmetrization}

Since the impedances of the faces of the scatterer are chosen to be equal, the problem can
be split into the symmetrical and antisymmetrical parts:
\begin{equation}
u^{\rm sc}(x,y) = u^{\rm a}(x,y) + u^{\rm s}(x,y),
\label{eq0201}
\end{equation}
where
\[
u^{\rm a}(x,y) = - u^{\rm a}(x,-y),
\qquad
u^{\rm s}(x,y) = u^{\rm s}(x,-y)
\]
are the antisymmetrical and symmetrical parts, respectively.

The symmetrical and antisymmetrical parts correspond to the incident waves
\[
u^{\rm in,s} = \frac{1}{2}[
\exp \{ - i k_0 (x \cos \theta^{\rm in} + y \sin \theta^{\rm in} ) \}
+
\exp \{ - i k_0 (x \cos \theta^{\rm in} - y \sin \theta^{\rm in} ) \}
],
\]
\[
u^{\rm in,a} = \frac{1}{2}[
\exp \{ - i k_0 (x \cos \theta^{\rm in} + y \sin \theta^{\rm in} ) \}
-
\exp \{ - i k_0 (x \cos \theta^{\rm in} - y \sin \theta^{\rm in} ) \}
],
\]
respectively.

The problems for $u^{\rm a}$ and $u^{\rm s}$ can be formulated as mixed boundary
value problems in the half-plane $y> 0$.
Boundary conditions for $u^{\rm a}$ are as follows:
\begin{equation}
\left[
\frac{\ptl}{\ptl y} - \eta
\right] u^{\rm a}(x,+0) = i k_0 \sin \theta^{\rm in} \exp \{ -i k_0 x \cos \theta^{\rm in} \}
\qquad   |x| < a,
\label{eq0202}
\end{equation}
\begin{equation}
u^{\rm a}(x,0) = 0,
\qquad   |x| > a.
\label{eq0203}
\end{equation}
Boundary conditions for $u^{\rm s}$ are as follows:
\begin{equation}
\left[
\frac{\ptl}{\ptl y} - \eta
\right] u^{\rm s}(x,+0) = \eta \exp \{ -i k_0 x \cos \theta^{\rm in} \}
\qquad   |x| < a,
\label{eq0204}
\end{equation}
\begin{equation}
\frac{\ptl}{\ptl y}u^{\rm s}(x,+0) = 0,
\qquad   |x| > a.
\label{eq0205}
\end{equation}

Below we study the symmetrical and the antisymmetrical problem separately (in parallel).
In both cases, we are interested in the field for $y > 0$ only.

The directivity of the scattered field is a sum of the symmetrical and antisymmetrical part:
\begin{equation}
S(\theta, \theta^{\rm in}) = S^{\rm s}(\theta, \theta^{\rm in}) + S^{\rm a}(\theta, \theta^{\rm in}),
\label{obvious01}
\end{equation}
where the last two values are defined similarly to (\ref{eq0106}).





\section{Local behavior of wave fields near the edges}

Here we study the growth of the solutions near the vertices. This growth is limited by the Meixner's conditions.

\begin{figure}[ht]
\centerline{\epsfig{file=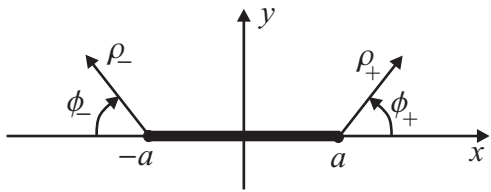}}
\caption{Local coordinates} \label{fig01}
\end{figure}

Introduce local cylindrical
variables $(\rho_\pm , \phi_\pm)$ (Fig.~\ref{fig01}).
Consider the {\em total\/} field in the {\bf antisymmetrical case}, i.\ e.\ consider the
function $u = u^{\rm a} + u^{\rm in, a}$.
 The Meixner's series
for a solution has form
\begin{equation}
u(\rho, \phi) = \sum_m \sum_n (k_0 \rho)^{\nu_m} \log^n (k_0 \rho) f_{m,n} (\phi),
\label{eq0301}
\end{equation}
where $\rho= \rho_{\pm}$, $\phi = \phi_{\pm}$, $f_{m,n} (\phi) = f^\pm_{m,n}(\phi_\pm)$.
This series is substituted into the Helmholtz equation and into the boundary conditions.
Also, some terms of the series are considered as prohibited according to the Meixner's condition mentioned above.
As the result, we get the following asymptotic expansion of the field:
\[
u = c (k_0 \rho)^{1/2} \sin(\phi/2) - \frac{2 c \eta}{3 \pi k_0 }(k_0 \rho)^{3/2} \phi \cos (3\phi/2)
\]
\begin{equation}
- \frac{2 c \eta}{3 \pi k_0 } (k_0 \rho)^{3/2} \log(k_0 \rho) \sin (3 \phi/2)
+ O(\log^2(k_0 \rho) (k_0 \rho)^{5/2}).
\label{eq0307}
\end{equation}
Now consider the {\bf symmetrical case}, i.\ e.\
let be  $u = u^{\rm s} + u^{\rm in, s}$.
The asymptotics for this case is as follows:
\begin{equation}
u = d - \frac{\eta d}{\pi} \rho \log (k_0 \rho) \cos(\phi) +
\frac{\eta d}{k_0 \pi}\rho \phi \sin(\phi) + O((k_0 \rho)^2 \log^2 (k_0 \rho) ).
\label{eq0305}
\end{equation}

Note that constants $c$ and $d$ in (\ref{eq0307}) and (\ref{eq0305})  are undetermined. Of course both constant take
different values for two edges, i.\ e.\ totally we introduce four constants $c_\pm$ and $d_\pm$ here.


\section{Formulation of Wiener--Hopf functional problems}

\subsection{Antisymmetrical case}

\begin{figure}[ht]
\centerline{\epsfig{file=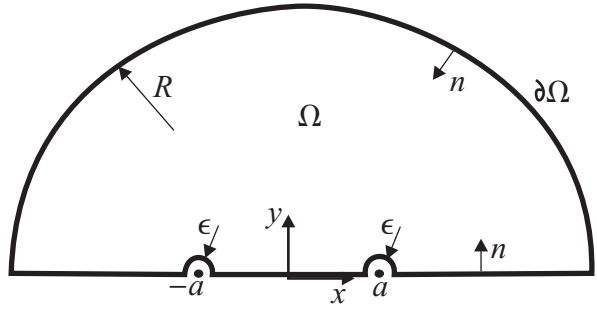}}
\caption{Contour for the Green's formula} \label{fig02}
\end{figure}

Consider domain $\Omega$ shown in Fig.~\ref{fig02}.
This domain is bounded by a part of $x$-axis, two small arcs (having radii $\epsilon \to 0$) encircling the vertices, and a large
arc (having radius $R \to \infty$) mimicking
the infinity. Consider two functions, both solutions of Hemholtz equation (\ref{eq0101}) in $\Omega$.
The first function
is $u^{\rm a}$ (the scattered field in the antisymmetrical case), and the second function is
an outgoing or decaying plane wave $w$:
\begin{equation}
w = w(k,x,y) =  \exp \left\{ i \left(k x + \xi(k) y \right)  \right\},
\label{eq0601}
\end{equation}
\begin{equation}
\xi(k) \equiv \sqrt{k_0^2 - k^2},
\label{eq0603b}
\end{equation}
where $k$ is a real value.
The branch of square root $\xi$ is chosen in such a way that while $|k| < {\rm Re}[k_0]$ the values of the square root are close to positive real. By continuity, the values of the square root for $|k| > {\rm Re}[k_0]$
are close to positive imaginary (the real axis
passes below the point~$k_0$ due to the limiting absorption principle). Note that $w$ is a solution of the Helmholtz equation for each
value of parameter $k$.

Apply the Green's formula to these two functions in $\Omega$:
\begin{equation}
\int_{\ptl \Omega} \left[
\frac{\ptl u^{\rm a}}{\ptl n} w - \frac{\ptl w}{\ptl n} u^{\rm a}
\right] dl = 0.
\label{eq0511}
\end{equation}
Since function $u^{\rm a}$ obeys the radiation condition, the integral over the large arc tends to zero as $R \to \infty$. The integrals over small arcs tend to zero as $\epsilon \to 0$ due to the local asymptotic expansions at the vertices.
Thus, only the integral over the
parts of the $x$-axis  should be considered.

Define the following values:
\begin{equation}
\check U_- (k)= \int \limits_{-\infty}^{-a} \left[
\frac{\ptl u^{\rm a}}{\ptl n} w - \frac{\ptl w}{\ptl n} u^{\rm a}
\right] dx =
\int \limits_{-\infty}^{-a}
\frac{\ptl u^{\rm a}(x,+0)}{\ptl y} e^{i k x} dx  ,
\label{eq0602}
\end{equation}
\begin{equation}
\check U_0 (k)= \int \limits_{-a}^a \left[
\frac{\ptl u^{\rm a}(x,+0)}{\ptl n} w(x,+0) - \frac{\ptl w (x,+0)}{\ptl n} u^{\rm a}(x,+0)
\right] dx ,
\label{eq0603}
\end{equation}
\begin{equation}
\check U_+ (k)= \int_a^{\infty} \left[
\frac{\ptl u^{\rm a}}{\ptl n} w - \frac{\ptl w}{\ptl n} u^{\rm a}
\right] dx =
\int \limits_a^{\infty}
\frac{\ptl u^{\rm a}(x,+0)}{\ptl y} e^{i k x} dx .
\label{eq0604}
\end{equation}
According to (\ref{eq0511}) the following functional equations are valid for all real $k$:
\begin{equation}
\check U_- (k) + \check U_0 (k) + \check U_+ (k) =0 .
\label{eq0605}
\end{equation}
Expression (\ref{eq0603}) can be transformed using (\ref{eq0202}):
\[
\check U_0(k) = (\eta - i \xi(k)) \int \limits_{-a}^a u^{\rm a}(x,+0) e^{i k x} dx +
\]
\begin{equation}
\frac{k_0 \sin \theta^{\rm in}}{k - k_*} \left(
\exp \{ i (k - k_*) a \}-
\exp \{ -i (k - k_*) a \}
\right),
\label{eq0603a}
\end{equation}
where
$$
k_* = k_0 \cos \theta^{\rm in}.
$$
Define the values
\begin{equation}
U_- (k) \equiv \check U_- (k) - \frac{k_0 \sin \theta^{\rm in}}{k - k_*}
\exp \{ -i (k - k_*) a \}
\label{eq0602d}
\end{equation}
\begin{equation}
U_0 (k) \equiv (\eta - i \xi(k)) \int \limits_{-a}^a u^{\rm a}(x,+0) e^{i k x} dx
\label{eq0603d}
\end{equation}
\begin{equation}
U_+ (k) \equiv \check U_+ (k) + \frac{k_0 \sin \theta^{\rm in}}{k - k_*}
\exp \{ i (k - k_*) a \}.
\label{eq0604d}
\end{equation}

According to (\ref{eq0605}) these values obey the functional equation
\begin{equation}
U_- (k) + U_0 (k) + U_+ (k) =0 .
\label{eq0605t}
\end{equation}

Functions $\check U_j$, $j = -,0,+$ are defined as Fourier transforms taken on some parts of the real axis.
Thus, standard theorems can be used to establish properties of these functions as well as the properties of $U_j$:

\begin{itemize}

\item[\bf Property 1 ]
Function $U_- (k)$ defined by (\ref{eq0602d}) and (\ref{eq0602})
can be analytically continued onto the whole
lower half-plane from the real axis, and it is regular there.
 Note that since we assume that $k_0$ has a negligibly small positive
imaginary part, the important point $k = - k_0$ belongs to the lower half-plane, and the function
$U_- (k)$  is regular at this point.

\item[\bf Property 2 ] Similarly, function $U_+ (k)$ defined by (\ref{eq0604d}) and
(\ref{eq0604}) can be analytically continued onto the whole upper half-plane including $k = k_0$, and it is
regular everywhere in the upper half-plane except a pole at $k = k_*$. At this pole function $U_+$ has
a prescribed residue equal to $k_0 \sin \theta^{\rm in}$.

\item[\bf Property 3 ] Function
\begin{equation}
\tilde U_0 (k) = \left( \eta - i \xi(k) \right)^{-1} U_0 (k)
\label{eq0808a}
\end{equation}
is regular on the whole complex plane $k$.

\item[\bf Property 4 ] Applying Watson's lemma to the integral representations (\ref{eq0602}), (\ref{eq0603}),
(\ref{eq0604}) we can get the following growth estimations as $|k| \to \infty$ in the domains of {\em a priori\/} regularity
of the unknown functions:
\begin{equation}
U_+(k) =  O(k^{-1/2} e^{i k a}),
\qquad {\rm Arg}[e^{- i \pi /2} k] \le \pi/2 ,
\label{eq0609a}
\end{equation}
\begin{equation}
U_-(k) = O(k^{-1/2} e^{-i k a}),
\qquad {\rm Arg}[e^{ i \pi /2} k] \le \pi/2 ,
\label{eq0610a}
\end{equation}
\begin{equation}
U_0(k) = O(k^{-1/2} e^{-i k a}),
\quad {\rm Arg}[e^{- i \pi /2} k] \le \pi/2 ,
\label{eq0611a}
\end{equation}
\begin{equation}
U_0(k) = O(k^{-1/2} e^{i k a}),
\qquad {\rm Arg}[e^{ i \pi /2} k] \le \pi/2 ,
\label{eq0612a}
\end{equation}
Note that estimations (\ref{eq0609a}), (\ref{eq0610a}) require some algebra to derive.

\end{itemize}

\begin{figure}[ht]
\centerline{\epsfig{file=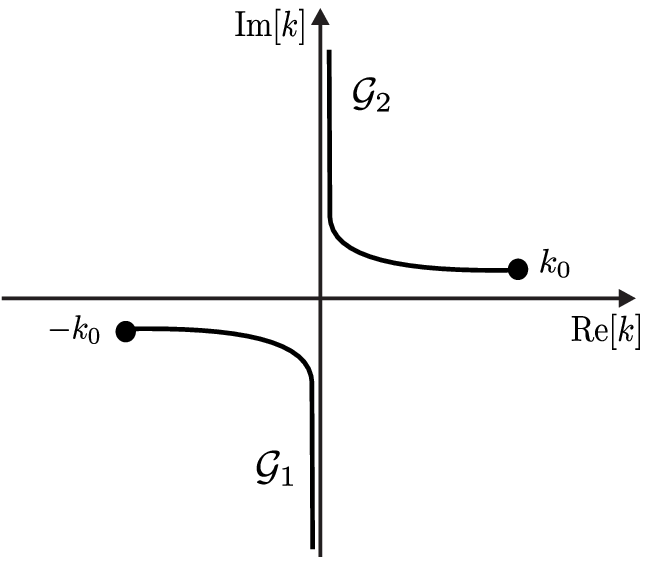}}
\caption{Cuts ${\cal G}_1$ and ${\cal G}_2$} \label{fig04_a}
\end{figure}

Introduce cuts ${\cal G}_1$ and ${\cal G}_2$ going from $-k_0$ and $k_0$ to infinity (see Fig.~\ref{fig04_a}).
These cuts go along the lines corresponding to the values of the square root $\pm \sqrt{k_0^2 - k^2}$ taken for real $k$.
Function $U_-$ can be naturally continued to the lower half-plane, function
$U_+$ can be naturally continued to the upper half-plane, and function $U_0$ can be continued to the
whole plane with the cuts ${\cal G}_1$ and ${\cal G}_2$. However, using  relations
\[
U_- (k) = - U_0 (k) - U_+ (k),
\qquad
U_+ (k) = - U_0 (k) - U_- (k)
\]
the functions $U_-$ can be continued to the upper half-plane with a cut ${\cal G}_2$, and
the functions $U_+$ can be continued to the lower half-plane with a cut ${\cal G}_1$. Moreover, it is
possible to study the Riemann surface of each function from the set $(U_-, U_+,U_0)$, and
prove that all branch points have order two and affixes~$\pm k_0$.

These properties enable us to formulate a functional problem for the functions $U_\pm$:

\begin{problem}
Find functions $U_+ (k)$, $U_- (k)$, regular in the complex plane
with the cuts ${\cal G}_1$ and ${\cal G}_2$, such that

\begin{itemize}

\item
function $U_-$ is regular in the lower half-plane;

\item
function $U_+$ is regular in the upper half-plane except a simple pole at $k = k_*$
with a residue equal to $k_0 \sin \theta^{\rm in}$;

\item
function  $(\eta - i \xi(k) )^{-1} U_0(k)$  is regular on the whole plane (here $U_0$ is defined as $U_0 \equiv - (U_+ + U_-)$);

\item
functions $U_+$, $U_-$, $\tilde U_0$ obey growth restrictions (\ref{eq0609a}), (\ref{eq0610a}),
(\ref{eq0611a}), (\ref{eq0612a}).

\end{itemize}

\label{functional_problem_A1}
\end{problem}

The formulation of the functional problem means that we forget about the definition of the unknown functions through the wave fields, and look for  functions $U_+ (k)$, $U_- (k)$ obeying
Problem~\ref{functional_problem_A1} and having arbitrary nature.

Let a solution of the functional problem be found. Let us describe the link between the directivity
$S^{\rm a} (\theta)$ for the antisymmetrical problem and the solution of the functional problem.
Apply Green's formula (\ref{eq0511}) to the domain $\Omega$, take
$u^{\rm a}$ as $u$, and
$u^{\rm in,a}(x,y)$ as $w$.
The integral over the large arc tends to a constant linked with the directivity.
The result is as follows:
\begin{equation}
S^{\rm a} (\theta,\theta^{\rm in}) = -e^{- i \pi /4} k_0 \sin \theta \,\,  \tilde U_0 (- k_0 \cos(\theta)).
\label{eq0617}
\end{equation}
Note that $\tilde U_0$ depends on $\theta^{\rm in}$ implicitly.


\subsection{Functional problem for the symmetrical case}

In the symmetrical case define functions $V_-(k)$, $V_+ (k)$, $V_0(k)$ by formulae
\begin{equation}
V_- (k) = \int_{-\infty}^{-a}
\exp \{ i k x \}  u^{\rm s} (x,+0) dx  - \frac{i}{k-k_*}\exp\{-i(k-k_*)a\},
\label{eq0602e}
\end{equation}
\begin{equation}
V_0 (k) =  \frac{ i \left(\eta - i \xi(k) \right)}{\eta  \xi(k)}
\int \limits_{-a}^a\exp \{ i k x \}\frac{\ptl{ u^{\rm s}(x, +0)}}{\ptl y}dx ,
\label{eq0603e}
\end{equation}
\begin{equation}
V_+ (k) =  \int_{a}^\infty
\exp \{ i k x \}
 u^{\rm s} (x,+0) dx  + \frac{i}{k-k_*}\exp\{i(k-k_*)a\},
\label{eq0604e}
\end{equation}
which are similar to (\ref{eq0602d}), (\ref{eq0603d}), (\ref{eq0604d}). A functional equation is valid
for these functions:
\begin{equation}
V_- (k) + V_0 (k) + V_+ (k) = 0.
\label{eq0605a}
\end{equation}
The growth estimations for the new unknown functions are as follows:
\begin{equation}
V_+(k) =  O(k^{-1} e^{i k a}),
\qquad {\rm Arg}[e^{- i \pi /2} k] \le \pi/2 ,
\label{eq0609b}
\end{equation}
\begin{equation}
V_-(k) = O(k^{-1} e^{-i k a}),
\qquad {\rm Arg}[e^{ i \pi /2} k] \le \pi/2 ,
\label{eq0610b}
\end{equation}
\begin{equation}
V_0(k) = O(k^{-1} e^{-i k a}),
\quad {\rm Arg}[e^{- i \pi /2} k] \le \pi/2 ,
\label{eq0611b}
\end{equation}
\begin{equation}
V_0(k) = O(k^{-1} e^{i k a}),
\qquad {\rm Arg}[e^{ i \pi /2} k] \le \pi/2 .
\label{eq0612b}
\end{equation}

The functional problem for the functions $V_\pm$ is as follows:

\begin{problem}
Find functions $V_+ (k)$, $V_+ (k)$, regular in the complex plane with the cuts ${\cal G}_1$ and ${\cal G}_2$, such that

\begin{itemize}

\item
function $ V_- (k)$ is regular in the lower half-plane;

\item
function $V_+ (k)$ is regular in the upper half-plane except a simple pole at $k = k_*$
with residue equal to $i$;

\item
function
\begin{equation}
\tilde V_0 =  \frac{\eta\xi(k)}{i(\eta - i \xi(k) )} V_0(k)
\end{equation}
 is regular in the whole plane (here $V_0$ is defined as $V_0 \equiv - (V_+ + V_-)$);

\item
functions $V_+$, $V_-$, $\tilde V_0$ obey growth restrictions (\ref{eq0609b}), (\ref{eq0610b}),
(\ref{eq0611b}), (\ref{eq0612b}).

\end{itemize}

\label{functional_problem_S1}
\end{problem}

The expression for the directivity of the symmetrical problem is as follows:
\begin{equation}
S^{\rm s} (\theta, \theta^{\rm in}) =  e^{- i \pi /4}  \tilde V_0 (- k_0 \cos(\theta)).
\label{eq0617a}
\end{equation}


\section{Auxiliary Wiener--Hopf functional problem and embedding formula}

\subsection{Auxiliary functions. Antisymmetrical problem}

Consider Problem~\ref{functional_problem_A1}. Here we modify this functional problem and formulate
a problem for the auxiliary functions. The following modifications are made. First, two pairs of
auxiliary functions are introduced. They are $(U^1_-, U^1_+)$, $(U^2_-, U^2_+)$.  This enables
us to construct a basis of solutions for a family of initial
functional problems indexed by parameter $\theta^{\rm in}$. Second, functions $U^{1,2}_+$ are required
to have no poles (i.\ e.\ the conditions of analyticity become more strict). Third, faster growth at infinity is allowed  (i.\ e.\ growth restriction become weaker).

\begin{problem}
Find functions $U^{1,2}_+ (k)$, $U^{1,2}_- (k)$, regular in the complex plane
with the cuts ${\cal G}_1$ and ${\cal G}_2$, such that

\begin{itemize}

\item
functions $U^{1,2}_-$ are regular in the lower half-plane;

\item
functions $U^{1,2}_+$ are regular in the upper half-plane;

\item
functions
\begin{equation}
\tilde U_0 = (\eta - i \xi(k) )^{-1} U^{1,2}_0(k)
\end{equation}
are regular on the whole plane  (here functions $U^{1,2}_0$ are defined as $U^{1,2}_0 \equiv - (U^{1,2}_+ + U^{1,2}_-)$);

\item
functions $U_+$, $U_-$, $\tilde U_0$ obey growth restrictions (\ref{eq0609c}), (\ref{eq0610c}),
(\ref{eq0611c}), (\ref{eq0612c}) formulated below.

\end{itemize}

\label{functional_problem_A2}
\end{problem}

The growth restrictions for this functional problem have the following form:
\begin{equation}
U^{j}_+(k) = \delta_{j,2} (e^{-i \pi /2 } k)^{1/2} e^{i k a} + O(k^{-1/2} e^{i k a}),
\qquad {\rm Arg}[e^{- i \pi /2} k] \le \pi/2 ,
\label{eq0609c}
\end{equation}
\begin{equation}
U^{j}_-(k) = \delta_{j,1} (e^{i \pi /2 } k)^{1/2} e^{-i k a} + O(k^{-1/2} e^{-i k a}),
\qquad {\rm Arg}[e^{ i \pi /2} k] \le \pi/2 ,
\label{eq0610c}
\end{equation}
\begin{equation}
\tilde U^{j}_0(k) = - \delta_{j,1} ( e^{- i \pi /2 } k)^{-1/2} e^{- i k a} + O(k^{-3/2} e^{-i k a}),
\quad {\rm Arg}[e^{- i \pi /2} k] \le \pi/2 ,
\label{eq0611c}
\end{equation}
\begin{equation}
\tilde U^{j}_0(k) = - \delta_{j,2} ( e^{i \pi /2 } k)^{-1/2} e^{i k a} + O(k^{-3/2} e^{i k a}),
\qquad {\rm Arg}[e^{ i \pi /2} k] \le \pi/2 ,
\label{eq0612c}
\end{equation}
where $j = 1,2$, and $\delta$ is the Kronecker's symbol.

Organize the solution of the auxiliary functional problem as a matrix
\begin{equation}
{\rm U}(k) = \left(
\begin{array}{cc}
U^1_-(k) & U^1_+ (k) \\
U^2_-(k) & U^2_+ (k)
\end{array}
\right).
\label{eq0613}
\end{equation}

Let us show that the solution of Problem~\ref{functional_problem_A2} is unique. Namely, let there exist two such solutions
${\rm U}$ and $\bar {\rm U}$. Consider the expression ${\rm J} = \bar {\rm U} {\rm U}^{-1}$.
This expression is equal to
\begin{equation}
{\rm J} = \frac{1}{D} \left(
\begin{array}{cc}
D_{1,1} & D_{1,2} \\
D_{2,1} & D_{2,2}
\end{array}
\right)
\label{eq0614}
\end{equation}
where
\[
D = | {\rm U} |,
\]
\[
D_{1,1} = \left|
\begin{array}{cc}
\bar U^1_-(k) & \bar U^1_+ (k) \\
U^2_-(k) & U^2_+ (k)
\end{array}
\right|,
\]
\[
D_{1,2} = \left|
\begin{array}{cc}
U^1_-(k) & U^1_+ (k) \\
\bar U^1_-(k) & \bar U^1_+ (k)
\end{array}
\right|,
\]
\[
D_{2,1} = \left|
\begin{array}{cc}
\bar U^2_-(k) & \bar U^2_+ (k) \\
U^2_-(k) & U^2_+ (k)
\end{array}
\right|,
\]
\[
D_{2,2} = \left|
\begin{array}{cc}
U^1_-(k) & U^1_+ (k) \\
\bar U^2_-(k) & \bar U^2_+ (k)
\end{array}
\right|,
\]
where $|\cdot|$ denotes determinant of the matrix.

All five determinants can be analyzed as follows. Consider $D$ as an example.
Study two representations of this determinant (they are equivalent due to linear dependence of
$U^j_-$, $U^j_+$, and $\tilde U^j_0$):
\begin{equation}
D = - ( \eta - i \xi(k) )
\left(
\begin{array}{cc}
U^1_- & \tilde U^1_0 \\
U^2_- & \tilde U^2_0
\end{array}
\right) =
-( \eta - i \xi(k) )
\left(
\begin{array}{cc}
\tilde U^1_0 & U^1_+  \\
\tilde U^2_0 & U^2_+
\end{array}
\right).
\label{eq0615}
\end{equation}
The first representation can be used to study the behaviour of
\[
\tilde D(k) \equiv - ( \eta - i \xi(k) )^{-1}  D(k)
\]
in the lower half-plane, and the second representation can be used to study the behaviour of the same
function in the upper half-plane. One can see that $\tilde D$ is analytical in both half-planes, and grows
as a constant equal to $-1$ in both half-planes. Thus,
according to Liouville's theorem,
\[
\tilde D \equiv -1.
\]
A similar reasoning can be applied to each of four other determinants. The result is
\[
{\rm J}(k) \equiv {\rm I},
\]
which is the identity matrix, i.\ e.\ the solution is unique.
Note that the determinant $D(k)$ can have no zeros except the zeros of the function $\eta - i \xi(k)$.

\subsection{Auxiliary functions. Symmetrical problem}

Similarly to the antisymmetrical case, introduce an auxiliary functional problem for the symmetrical case.

\begin{problem}
Find functions $V^{1}_+ (k)$, $V^{2}_+ (k)$, $V^{1}_- (k)$, $V^{2}_- (k)$, regular in the complex plane
with the cuts ${\cal G}_1$ and ${\cal G}_2$, such that

\begin{itemize}

\item
functions $V^{j}_-$ are regular in the lower half-plane;

\item
functions $V^{j}_+$ are regular in the upper half-plane;

\item
functions
\begin{equation}
\label{eq0605_s}
\tilde V^{j}_0 \equiv -\frac{  \xi(k)}{ i \left(\eta - \i \xi(k) \right)} (V^{j}_- + V^{j}_+)
\end{equation}
are regular on the whole plane;

\item
functions $V^{j}_+$, $V^{j}_-$, $\tilde V^{j}_0$ obey growth restrictions (\ref{eq0609s}), (\ref{eq0610s}), (\ref{eq0611s}), (\ref{eq0612s})  formulated below.

\end{itemize}

\label{functional_problem_S2}
\end{problem}
The growth conditions for this functional problem have the following form:
\begin{equation}
V^{j}_+(k) = \delta_{j,2}  e^{i k a} + O(k^{-1}e^{i k a}),
\qquad {\rm Arg}[e^{- i \pi /2} k] \le \pi/2 ,
\label{eq0609s}
\end{equation}
\begin{equation}
V^{j}_-(k) = \delta_{j,1}  e^{-i k a} + O(k^{-1} e^{-i k a}),
\qquad {\rm Arg}[e^{ i \pi /2} k] \le \pi/2 ,
\label{eq0610s}
\end{equation}
\begin{equation}
\tilde V^{j}_0(k) = - \delta_{j,1} e^{- i k a} + O(k^{-1} e^{-i k a}),
\quad {\rm Arg}[e^{- i \pi /2} k] \le \pi/2 ,
\label{eq0611s}
\end{equation}
\begin{equation}
\tilde V^{j}_0(k) = - \delta_{j,2} e^{i k a} + O(k^{-1}e^{i k a}),
\qquad {\rm Arg}[e^{ i \pi /2} k] \le \pi/2 .
\label{eq0612s}
\end{equation}
The solution of the functional problem can be organized as a matrix
\begin{equation}
{\rm V}(k) = \left(
\begin{array}{cc}
V^{1}_-(k) & V^{1}_+ (k) \\
V^{2}_-(k) & V^{2}_+ (k)
\end{array}
\right).
\label{eq0613_s}
\end{equation}

Using representation similar to (\ref{eq0614}) one can show that
Problem~\ref{functional_problem_S2} has a unique solution.


\subsection{Embedding formula}

Consider the {\bf antisymmetrical} case.
Let row vector $(U_- , U_+)$ be a solution of Problem~\ref{functional_problem_A1}, and let ${\rm U}(k)$ be
a solution of Problem~\ref{functional_problem_A2} in the matrix form (\ref{eq0613}).
Find functions $r_1(k)$ and $r_2(k)$ such that
\begin{equation}
(U_-(k), U_+(k)) = (r_1(k) , r_2(k))
\left(
\begin{array}{cc}
U^1_-(k) & U^1_+ (k) \\
U^2_-(k) & U^2_+ (k)
\end{array}
\right).
\label{eq2001}
\end{equation}
Due to Cramer's rule,
\begin{equation}
r_1 = \frac{D_1}{D}, \qquad r_2 = \frac{D_2}{D},
\label{eq2002}
\end{equation}
where
\begin{equation}
D_1 =
\left|
\begin{array}{cc}
 U_-(k) &  U_+ (k) \\
U^2_-(k) & U^2_+ (k)
\end{array}
\right|,
\qquad
D_2 =
\left|
\begin{array}{cc}
 U^1_-(k) &  U^2_+ (k) \\
U_-(k) & U_+ (k)
\end{array}
\right|.
\label{eq2003}
\end{equation}

Determinant $D$ was calculated in the previous section using representation (\ref{eq0615}).
Determinants $D_1$, $D_2$ can be analyzed similarly to determinant $D$, namely there
exist two representations for each determinant enabling one to study these determinants in the upper and lower
half-plane:
\begin{equation}
D_1 = - ( \eta - i \xi(k) )
\left(
\begin{array}{cc}
U_- & \tilde U_0 \\
U^2_- & \tilde U^2_0
\end{array}
\right) =
-( \eta - i \xi(k) )
\left(
\begin{array}{cc}
\tilde U_0 & U_+  \\
\tilde U^2_0 & U^2_+
\end{array}
\right),
\label{eq2003add01}
\end{equation}
\begin{equation}
D_2 = - ( \eta - i \xi(k) )
\left(
\begin{array}{cc}
U^1_- & \tilde U^1_0 \\
U_- & \tilde U_0
\end{array}
\right) =
-( \eta - i \xi(k) )
\left(
\begin{array}{cc}
\tilde U^1_0 & U^1_+  \\
\tilde U_0 & U_+
\end{array}
\right).
\label{eq2003add02}
\end{equation}
 Using these representations and applying the Liouville's theorem one can prove that
\begin{equation}
D_1 = \frac{\left( \eta - i \sqrt{k_0^2 - k^2} \right)}{k - k_*}R_1,
\label{eq2003add03}
\end{equation}
\begin{equation}
D_2 = \frac{\left( \eta - i \sqrt{k_0^2 - k^2} \right)}{k - k_*}R_2,
\label{eq2003add04}
\end{equation}
where  $R_1$, $R_2$ are some constants. $R_1$, $R_2$ can be obtained by calculating residues of determinants $D_1$, $D_2$ at
the point $k=k_*$. These residues can be found either from (\ref{eq2003add01}), (\ref{eq2003add02}) or from
(\ref{eq2003add03}), (\ref{eq2003add04}). Comparing these representations, obtain
\begin{equation}
R_1 = -\sqrt{k_0^2 - k_*^2}\, \tilde U^2_0(k_*), \quad R_2 = \sqrt{k_0^2 - k_*^2} \, \tilde U^1_0(k_*).
\end{equation}
Substituting $r_1$ and $r_2$ into (\ref{eq2003}) obtain the embedding formula:
\begin{equation}
\tilde U_{0}(k,k_*) =\frac{\xi(k_*)}{k-k_*}\left(\tilde U_0^{1}(k_*)\tilde U^{2}_0(k)  - \tilde U_0^{1}(k)\tilde U^{2}_0(k_*)\right).
\label{embedding_a}
\end{equation}
According to embedding formula we can focus our efforts on finding the solution of
Problem~\ref{functional_problem_A2}, namely on  functions  $U^{j}_0(k)$, $j= 1,2$.

Conducting a similar procedure one can obtain an embedding formula for the {\bf symmetrical case}:
\begin{equation}
\tilde V_{0}(k,k_*) =\frac{i\eta}{(k-k_*)}\left(\tilde V^{2}_0(k_*)\tilde V_0^{1}(k) - \tilde V^{2}_0(k)\tilde V_0^{1}(k_*)  \right).
\label{embedding_s}
\end{equation}


\section{Matrix Riemann--Hilbert  formulation for auxiliary functional problems}
\subsection{Antisymmetrical problem}

Here we present a matrix Riemann--Hilbert formulation for the antisymmetrical case.

\begin{figure}[ht]
\centerline{\epsfig{file=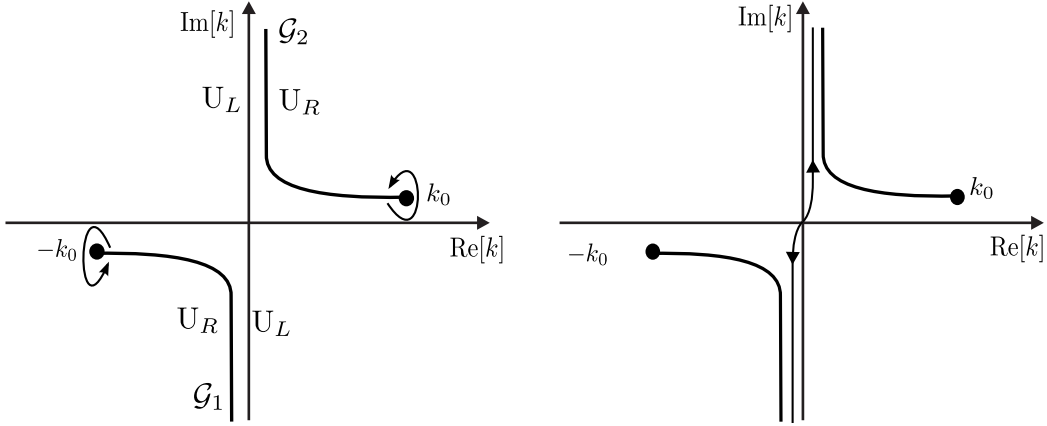}}
\caption{(left) Bypasses around $k_0$ and $-k_0$. (right) Analytical continuation of
the square roots}
\label{fig04}
\end{figure}

Let us make some preliminary steps. Consider the cuts ${\cal G}_1$ and ${\cal G}_2$  (see Fig.~\ref{fig04}, left). The values on the left shores (when going from $\pm k_0$ to $\infty$) of the cuts are denoted by symbols with lower index $L$;
the values on the right shores are denoted by index~$R$.

Consider the bypasses about $\pm k_0$ and  going from a point on the left shore to the right shore,
i.\ e.\ going in the positive direction. Our current aim is to describe the transformation of the
matrix ${\rm U}$ occuring as a result of  the bypass.
Namely, let us prove that
\begin{equation}
{\rm U}_R(k) = {\rm U}_L(k)\, {\rm M}_1(k), \qquad k \in {\cal G}_1 ,
\label{eq0801}
\end{equation}
\begin{equation}
{\rm U}_R(k) = {\rm U}_L(k)\, {\rm M}_2(k), \qquad k \in {\cal G}_2 ,
\label{eq0802}
\end{equation}
with
\begin{equation}
{\rm M}_1 (k) = \left(  \begin{array}{cc}
1 & 2 i \xi/(\eta - i \xi) \\
0 & (\eta + i \xi)/(\eta - i \xi)
\end{array} \right) ,
\label{eq0806}
\end{equation}
\begin{equation}
{\rm M}_2 (k) = \left(  \begin{array}{cc}
(\eta + i \xi)/(\eta - i \xi) & 0 \\
2 i \xi/(\eta - i \xi)      & 1
\end{array} \right) .
\label{eq0805}
\end{equation}

The analytic continuation of the square root
$\xi(k)\equiv \sqrt{k_0^2 - k^2}$ on the cuts ${\cal G}_{1,2}$ is defined as follows. This square root is equal to $k_0$ for $k =0$.
Then, introduce the paths shown in Fig.~\ref{fig04} (right). These paths go from zero to the left shores
of ${\cal G}_{1,2}$. The values of the square root on ${\cal G}_{1,2}$
is taken as the result of the continuation along these paths. The values of the square root are taken for
${\rm M}_{1,2}$ from the left shores.

Derive (\ref{eq0802}).
Consider contour ${\cal G}_{2}$ associated with matrix ${\rm M}_2$.
Continue functional equation (\ref{eq0605}):
\begin{equation}
(U_-^j(k))_L = - U^j_+(k) - \left(\eta - i \xi(k) \right)\tilde U^j_0 (k),
\label{eq0803}
\end{equation}
\begin{equation}
(U_-^j(k))_R = - U^j_+(k) - \left(\eta + i \xi(k)\right)\tilde U^j_0 (k).
\label{eq0804}
\end{equation}
Then,
\[
(U_-^j (k))_R = \frac{\eta + i \xi(k)}{\eta - i \xi(k)} (U_-^j(k))_L +
\frac{2 i \xi(k)}{\eta - i \xi(k)} U_+^j (k).
\]
Note that functions $U^j_+$ and $\tilde U_0^j$ are not labeled as $R$ or $L$, since they do not change their
values after the considered bypass.
Thus, relations (\ref{eq0802}) and (\ref{eq0805}) are valid. Similarly one can prove
(\ref{eq0801}) and (\ref{eq0806}).

 Reformulate the
growth restrictions (\ref{eq0611c}) and (\ref{eq0612c}) according to (\ref{eq0605}) as follows:
\begin{equation}
U^{j}_- = i \, \delta_{j,1} ( e^{- i \pi /2 } k)^{1/2} e^{- i k a} + O(k^{-1/2} e^{-i k a}),
\quad {\rm Arg}[e^{- i \pi /2} k] \le \pi/2 ,
\label{eq0807}
\end{equation}
\begin{equation}
U^{j}_+ = i \, \delta_{j,2} ( e^{  i \pi /2 } k)^{1/2} e^{i k a} + O(k^{-1/2} e^{i k a}),
\qquad {\rm Arg}[e^{ i \pi /2} k] \le \pi/2 .
\label{eq0808}
\end{equation}
Both restrictions are related to the continuations along the paths shown in Fig.~\ref{fig04}.

Now we can formulate a Riemann--Hilbert problem for ${\rm U}$:

\begin{problem}
Find a matrix function ${\rm U}(k)$ of elements (\ref{eq0613}) such that

\begin{itemize}

\item
it is regular on the plane cut along the lines ${\cal G}_{1,2}$;

\item
it obeys functional equations (\ref{eq0801}), (\ref{eq0802}) with coefficients
(\ref{eq0805}), (\ref{eq0806}) on the cuts;

\item
it obeys growth restrictions (\ref{eq0609c}), (\ref{eq0610c}), (\ref{eq0807}), (\ref{eq0808});

\item
functions $U^j_+ (k) + U^j_- (k)$, $j = 1,2$ have zeros at
$k = k' \equiv \sqrt{k_0^2 + \eta^2}$;

\item
functions $U^j_{\pm}$ grow no faster than a constant near the points $\pm k_0$.

\end{itemize}

\label{WHH_with_zeros}
\end{problem}

The fourth condition (concerning zeros at $\pm k'$) are difficult to take into account,
so we would like to eliminate it.
Consider Riemann surface of the function $\sqrt{k_0^2 - k^2}$ cut along the lines
${\cal G}_{1,2}$. The surface is split into two sheets by the cuts. The sheet to which the point
$\sqrt{k_0^2 - 0^2} = k_0$ belongs will be called the physical sheet.
Consider the function $\eta - i \sqrt{k_0^2 - k^2}$ on this surface. Note that
this function has two zeros only on one sheet (on the physical one or on the other one).
If the zeros belong to the physical sheet, deform the contours ${\cal G}_{1,2}$ such that:

\begin{itemize}

\item
the end points remain the same;

\item
contour ${\cal G}_2$ remains symmetrical to ${\cal G}_1$ with respect to zero;

\item
zeros of $\eta - i \sqrt{k_0^2 - k^2}$ finally become not belonging to the physical sheet.

\end{itemize}

\begin{figure}[ht]
\centerline{\epsfig{file=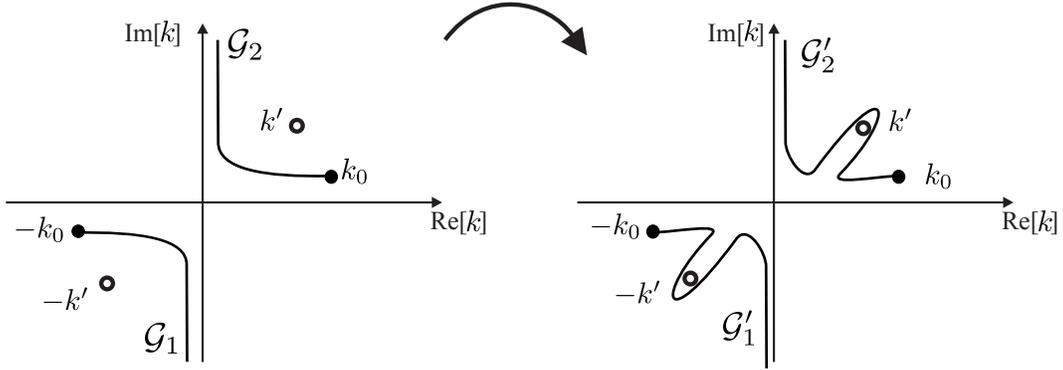,width = 14cm}}
\caption{Deformation of the cuts ${\cal G}_{1,2}$}
\label{fig05}
\end{figure}

A scheme of such contour deformation is shown in Fig.~\ref{fig05}.

If the zeros do not belong to the physical sheet from the very beginning, then no deformation
is needed. The domain of $\eta$ for which the zeros of
$\eta - i \sqrt{k_0^2 - k^2}$
belong to the physical sheet (and the deformation is needed) is
\begin{equation}
{\rm Im}[\eta] <0,\quad {\rm Re}[\eta] <0,
\end{equation}
i.\ e.\ it is the third quadrant of the complex plane.

Denote the resulting contours (deformed if the deformation is needed or undeformed otherwise)
by ${\cal G}_{1,2}'$.

{\bf Remark. }
Positions of the points $k'$ on the Riemann surface of $\sqrt{k_0^2 - k^2}$
can be found from condition (\ref{eq0104}). Namely, the boundary between the allowed values of $\eta$
and prohibited values is the real axis. Consider the function $k' = k' (\eta)$. This function
maps the real axis of $\eta$ into the parts ${\cal G}_1'' = (-\infty , - k_0 )$, ${\cal G}_2''= (k_0, \infty)$
of the real axis.
Consider the Riemann surface of $\sqrt{k_0^2 - k^2}$ cut along ${\cal G}_{1,2}''$. The surface will be split into
two sheets. Again, call the sheet containing the point $\sqrt{k_0^2 -0^2} = k_0$ the physical sheet.
The boundary ${\rm Im}[\eta] = 0$ corresponds to the cuts ${\cal G}_{1,2}''$.
The area  ${\rm Im}[\eta] < 0$ corresponds to the unphysical sheet.


Formulate the functional problem for the contours ${\cal G}_{1,2}'$.
According to the principles of analytical continuation, relations (\ref{eq0801}), (\ref{eq0802})
remain valid with the same matrices (\ref{eq0805}), (\ref{eq0806}). Thus, the formulation of the problem is almost the same:

\begin{problem}
Find a matrix function ${\rm U}(k)$ of elements (\ref{eq0613}) such that

\begin{itemize}

\item
it is regular on the plane cut along the lines ${\cal G}_{1,2}'$;

\item
it obeys functional equations (\ref{eq0801}), (\ref{eq0802}) with coefficients
(\ref{eq0805}), (\ref{eq0806}) on the cuts;

\item
it obeys growth restrictions (\ref{eq0609c}), (\ref{eq0610c}), (\ref{eq0807}), (\ref{eq0808});

\item
functions $U^j_{\pm}$ grow no faster than a constant near the points $\pm k_0$.

\end{itemize}

\label{WHH}
\end{problem}

\subsection{Symmetrical problem}

Similarly to the antisymmetrical case, there are two functional equations describing the transformation
of unknown functions at the cuts:
\begin{equation}
{\rm  V}_R(k) = {\rm  V}_L {\rm N}_1(k), \qquad k \in {\cal G}_1 ,
\label{eq0801sym}
\end{equation}
\begin{equation}
{\rm  V}_R(k) = {\rm  V}_L {\rm N}_2(k), \qquad k \in {\cal G}_2 .
\label{eq0802sym}
\end{equation}
 \begin{equation}
{\rm  N}_1 (k) = \left(  \begin{array}{cc}
1 & -2 \eta/(\eta - i \zeta) \\
0 & (\eta + i \zeta)/(i \zeta - \eta  )
\end{array} \right) ,
\label{eq0806sym}
\end{equation}
\begin{equation}
{\rm N}_2 (k) = \left(  \begin{array}{cc}
(\eta + i \zeta)/(i \zeta - \eta) & 0 \\
-2 \eta/(\eta - i \zeta)      & 1
\end{array} \right) ,
\label{eq0805sym}
\end{equation}

Reformulate growth restrictions (\ref{eq0611s}), (\ref{eq0612s}) according to (\ref{eq0605_s}) as follows:
\begin{equation}
 V^{j}_- =  \, \delta_{j,1}e^{- i k a} + O(k^{-1} \log(k)e^{-i k a}),
\quad {\rm Arg}[e^{- i \pi /2} k] \le \pi/2  ,
\label{eq0807sym}
\end{equation}
\begin{equation}
 V^{j}_+ =  \, \delta_{j,2} e^{i k a} + O(k^{-1} \log(k) e^{i k a}),
\quad {\rm Arg}[e^{ i \pi /2} k] \le \pi/2  .
\label{eq0808sym}
\end{equation}

Finally, formulate a functional problem for ${\rm V}$.

\begin{problem}
Find a matrix function ${\rm  V}(k)$ of elements (\ref{eq0613_s}) such that

\begin{itemize}

\item
it is regular on the plane cut along the lines ${\cal G}_{1,2}$;

\item
it obeys functional equations (\ref{eq0801sym}), (\ref{eq0802sym}) with coefficients
(\ref{eq0806sym}), (\ref{eq0805sym}) on the cuts;

\item
it obeys growth restrictions (\ref{eq0609s}), (\ref{eq0610s}), (\ref{eq0807sym}), (\ref{eq0808sym});

\item
functions $ V^{j}_{\pm}$ grow no faster than $(\sqrt{k_0 \mp k})^{-1/2}$ near the points $\pm k_0$.

\end{itemize}

\label{WHH_sym}
\end{problem}

\section{Conclusion}
The problem of diffraction by impedance strip is symmetrized and reduced to two Wiener--Hopf
functional problems
(Problem \ref{functional_problem_A1} and \ref{functional_problem_S1})
leading to directivities $S^{\rm a}(\theta,\theta^{\rm in})$ and
$S^{\rm s}(\theta,\theta^{\rm in})$.
Then auxiliary functional problems (Problem~\ref{functional_problem_A2} and~\ref{functional_problem_S2}) are introduced.
Using embedding formulae (\ref{embedding_a}) and (\ref{embedding_s}) a
simple connection with  Problem \ref{functional_problem_A1} and \ref{functional_problem_S1} is established.
Riemann--Hilbert problems (Problem~\ref{WHH} and \ref{WHH_sym}) for auxiliary solutions are formulated.

In the second part of the paper the family of Riemann--Hilbert problems indexed by an artificial parameter
will be introduced. A differential equation will be built with respect to this parameter. A novel technique of
OE--equation will be applied to solve this equation and find the solution of original problem.
Some numerical results will be presented.

\section*{Acknowledgements}

The work is supported by the grants RFBR 14-02-00573, Scientific Schools-283.2014.2,
RF Government grant 11.G34.31.0066.

The authors are grateful to participants of the seminar on wave diffraction held in S.Pb. branch of Steklov
Mathematical Institute of RAS (the chairman is Prof. V.\ M.\ Babich) for interesting discussions.

\end{document}